\documentclass{article}

\usepackage{amssymb}

\newcommand{\proof}{{\bf Proof:  }}
\newcommand{\remark}{{\bf Remark:  }}

\newcommand{\example}{{\bf Example:  }}

\newcommand{\dimv}{\underline{\dim}}

\newcommand{\hb}{\newline\hspace*{\fill}$\Box$}

\newtheorem{theorem}{Theorem}[section]
\newtheorem{lemma}[theorem]{Lemma}
\newtheorem{definition}[theorem]{Definition}
\newtheorem{proposition}[theorem]{Proposition}
\newtheorem{corollary}[theorem]{Corollary}

\begin{document}

\parindent0pt

\title{\bf Counting rational points of quiver moduli}

\author{Markus Reineke\\ Mathematisches Institut\\ Universit\"at M\"unster\\ D - 48149 M\"unster, Germany\\ e-mail: reinekem@math.uni-muenster.de}

\date{}

\maketitle

\begin{abstract} It is shown that rational points over finite fields of moduli spaces of stable quiver representations are counted by polynomials with integer coefficients. These polynomials are constructed recursively using an identity in the Hall algebra of a quiver.
\end{abstract}

\section{Introduction}\label{introduction}

Geometric Invariant Theory \cite{M} provides a general technique for formulation of interesting moduli problems for objects in an abelian category. The essence of this approach is to restrict the class of objects to be parametrized up to isomorphism to certain stable ones. The resulting moduli spaces are seldomly proper, making the determination of topological invariants a very difficult problem.\\[1ex]
Moduli spaces for stable representations of quivers, introduced in \cite{Ki}, provide an interesting testing ground for techniques of moduli theory, since they are easy to construct. They are of course also interesting in themselves, particularly in view of the classification problem of quiver representations. While these moduli spaces are always smooth, they are projective only in the case of indivisible dimension types for quivers without oriented cycles (in which case their Betti numbers were determined in \cite{HNS}).\\[1ex]
In this paper, a first step towards the determination of geometric invariants of quiver moduli in the general case is made. The main result asserts that rational points over finite fields of arbitrary quiver moduli are counted by polynomials with integer coefficients (Theorem \ref{polynomiality}). These polynomials are recursively computable (Theorem \ref{doublerecursion}), and their evaluation at one equals the Euler characteristic of the moduli space.\\[1ex]
In more algebraic terms, these polynomials count isomorphism classes of absolutely stable quiver representations over finite fields (to be discussed in section \ref{arithmetics}). Note that, as a particular example, one can count absolutely simple representations of free algebras in this way (Theorem \ref{freealgebra}). Although the recursive description Theorem \ref{doublerecursion} of the counting polynomials is quite complicated, it nevertheless opens the opportunity for extensive experimentation (see section \ref{examples} for some examples).\\[2ex]
Surprisingly, there seems to be no way to establish the recursive formula Theorem \ref{doublerecursion} by a direct counting procedure. Instead, it is derived (in section \ref{key}) as a consequence of an identity in (a completion of) the Hall algebra \cite{RH} of the quiver (section \ref{identity}). This proves the existence of rational functions counting rational points over finite fields; their polynomiality is established using basic arithmetic geometry (section \ref{counting}).\\[3ex]
{\bf Acknowledgments:} The author wants to thank M. Van den Bergh for pointing out certain techniques of \cite{Kac}, and K. Bongartz, P. Caldero, C. Deninger, A. Hubery, B. Keller, L. Le Bruyn and A. Schofield for useful discussions. 

\section{Quiver moduli}\label{quivermoduli}

Let $Q$ be a quiver with finite set of vertices $I$ and finitely many arrows $\alpha:i\rightarrow j$ between each pair of vertices $i,j\in I$. An element $d=\sum_{i\in I}d_ii\in{\bf N}I$ of the free abelian semigroup generated by $I$ will be called a dimension type. We introduce a bilinear form on ${\bf Z}I$, called Euler form, by $\langle d,e\rangle=\sum_{i\in I}d_ie_i-\sum_{i\rightarrow j}d_ie_j$.\\[1ex]
A linear function $\Theta:{\bf Z}I\rightarrow{\bf Z}$ will be called a stability. Associated to $\Theta$, we have a slope function $\mu:{\bf N}I\setminus 0\rightarrow{\bf Q}$ defined by $\mu(d):=\frac{\Theta(d)}{\dim d}$, where $\dim d=\sum_{i\in I}d_i$ denotes the total dimension of $d$. For a rational number $\mu\in{\bf Q}$, we denote by ${\bf N}I_\mu\subset{\bf N}I$ the set of all dimension types $d$ of slope $\mu(d)=\mu$. Note that ${\bf N}I_\mu$ is a finitely generated subsemigroup of ${\bf N}I$.\\[1ex]
For $k$ an arbitrary field, we denote by $kQ$ the path algebra of $Q$ over $k$. Then the abelian $k$-linear category of $k$-representations of $Q$ is equivalent to the category ${\rm mod} kQ$ of representations of $kQ$ (for all general notions and facts on representations of quivers the reader is referred to \cite{ARS,RiB}). The symbols ${\rm Hom}$, ${\rm End}$, ${\rm Aut}$ are understood as homomorphisms (resp. endomorphisms, automorphisms) in the category ${\rm mod} kQ$.\\[1ex]
A representation $M$ of $kQ$ of slope $\mu(M)=\mu(\dimv M)=\mu\in{\bf Q}$ is called semistable (resp.~stable) if for all non-zero proper subrepresentations $U$ of $M$, we have $\mu(U)\leq\mu(M)$ (resp.~$\mu(U)<\mu(M)$). It is called 
polystable of slope $\mu$ if it is isomorphic to a direct sum of stable representations of slope $\mu$. We list some basic facts on these notions from e.g. \cite{HNS}:\\[1ex]
The full subcategory ${\rm mod}_\mu kQ$ of semistable representations of $kQ$ of slope $\mu$ forms an exact abelian subcategory of ${\rm mod} kQ$, that is, it is closed under kernels, images, cokernels, and extensions. The simple objects in ${{\rm mod}}_\mu kQ$ are precisely the stable representations, thus the semisimple objects in ${{\rm mod}} kQ$ are precisely the polystable representations. In particular, we have the notion of the ${{\rm mod}}_\mu kQ$-relative socle of a semistable representation with all its usual properties.\\[1ex]
Any representation $M$ of $kQ$ possesses a unique Harder-Narasimhan filtration, that is, a filtration $0=M_0\subset M_1\subset\ldots\subset M_s=M$ such that each subquotient $M_i/M_{i-1}$ is semistable of some slope $\mu_i$, and $\mu_1>\ldots>\mu_s$.\\[2ex]
The $k$-representations of $Q$ of dimension type $d$ are parametrized by the affine space $R_d(k):=\bigoplus_{\alpha:i\rightarrow j}{\rm Hom}_k(k^{d_i},k^{d_j})$, on which the linear algebraic group $G_d(k):=\prod_{i\in I}{\rm GL}_{d_i}(k)$ acts via base change
$$(g_i)_i\cdot(M_\alpha)_\alpha:=(g_jM_\alpha g_i^{-1})_{\alpha:i\rightarrow j}.$$
Setting $\mu=\mu(d)$, we denote the set of points of $R_d$ corresponding to semistable (resp.~stable) representations of slope $\mu$ by $R^{ss}_d(k)$ (resp.~$R_d^{s}(k)$). We have open inclusions $R_d^{s}(k)\subset R_d^{ss}(k)\subset R_d(k)$.\\[1ex]
Now let $k$ be algebraically closed. It is proved in \cite{Ki} that $R_d^{s}$ is precisely the set of Mumford-stable points in $R_d(k)$ (for a certain $G_d(k)$-linearization of the trivial line bundle), so that there exists a geometric quotient $M^s_d(Q)(k):=R_d^{s}(k)/G_d(k)$. Moreover, the action of the quotient $PG_d(k)$ of $G_d(k)$ by the scalar matrices on $R_d^{s}(k)$ is scheme-theoretically free by \cite[Lemma 6.5]{HNS}, so that the quotient map turns $R_d^{s}(k)$ into a principal $PG_d(k)$-bundle by \cite[0.9]{M}.

\section{An identity in the Hall algebra}\label{identity}

We consider a completed version of the Hall algebra \cite{RH} of the quiver $Q$ and a certain map, called the integral, to a $q$-commutative formal power series ring. These will function as the main organizing tool for the counting procedure. The principle is the following: we will construct an identity between generating functions in the Hall algebra, which, via the above-mentioned integral, will yield a counting formula inside the $q$-commutative formal power series ring.\\[1ex]
The following defines a completed version of the usual Hall algebra of a quiver.

\begin{definition} Let $k$ be a finite field with $q$ elements, and let $Q$ be a finite quiver. Define $H((kQ))$ as the direct product of one-dimensional ${\bf Q}$-vector spaces with basis elements $[M]$, indexed by the isomorphism classes of finite dimensional representations of $kQ$, and with the following multiplication: $[M]\cdot[N]=\sum_{[X]}F_{M,N}^X\cdot[X]$, where $F_{M,N}^X$ denotes the number of $kQ$-subrepresentations $U$ of $X$ such that $U\simeq N$ and $X/U\simeq M$.
\end{definition}

This definition is standard except for the use of the direct product. To see that the multiplication is well-defined, it suffices to observe that representations $M$ and $N$ in whose product a fixed representation $X$ appears with non-zero coefficient already have smaller dimension than $X$, and there are only finitely many such since we work over a finite field $k$.\\[1ex]
Note furthermore that the direct product $\prod_{[M]\in{\rm mod}_\mu kQ}{\bf Q}\cdot[M]$ over all isomorphism classes in ${\rm mod}_\mu kQ$ defines a subalgebra of $H((kQ))$.

\begin{definition} Define ${\bf Q}_q[[I]]$ as the direct product of the one-dimensional ${\bf Q}$-vector spaces with basis elements $t^d$ indexed by $d\in{\bf N}I$, and with multiplication $t^d\cdot t^e=q^{-\langle d,e\rangle}t^{d+e}$.
\end{definition}

The following is \cite[Lemma 6.3]{HNS}.

\begin{lemma}\label{propint} The map sending a basis element $[M]$ to the element $$\int[M]:=\frac{1}{|{\rm Aut}(M)|}t^{\dimv M}$$ induces a ${\bf Q}$-algebra homomorphism $\int:H((kQ))\rightarrow {\bf Q}_q[[I]]$.
\end{lemma}

\proof By a formula of C. Riedtmann \cite{Rie}, we 
have
$$F_{M,N}^X=q^{-\dim{\rm Hom}(M,N)}\frac{|{\rm Aut}(X)|}{|{\rm 
Aut}(M)|\cdot |{\rm Aut}(N)|}\cdot|{\rm Ext}^1(M,N)_X|,$$
where ${\rm Ext}^1(N,M)_X$ denotes the set of extension classes corresponding to short exact sequences with middle 
term isomorphic to $X$. The lemma now follows from a direct calculation, using the definition of ${\bf Q}_q[[I]]$ and the identity $\langle\dimv M,\dimv N\rangle=\dim{\rm Hom}(M,N)-\dim{\rm Ext}^1(M,N)$.\hb

For $\mu\in{\bf Q}$, define an element $e_\mu$ of $H((kQ))$ by $$e_\mu=\sum_{[M]\in{\rm mod}_\mu kQ}[M].$$
Denote by ${\cal S}_\mu$ the set of isomorphism classes $[S]$ of stable representations $S$ in ${\rm mod}_\mu kQ$.

\begin{lemma}\label{keyidentity} Writing $e_\mu^{-1}=\sum_{[M]\in{\rm mod}_\mu kQ}\gamma_M[M]$, we have $\gamma_M=0$ if $M$ is not po\-ly\-stable of slope $\mu$, and
$$\gamma_M=\prod_{[S]\in{\cal S}_\mu}(-1)^{m_S}|{\rm End}(S)|^{m_S\choose 2}$$
if $M=\bigoplus_{[S]\in{\cal S}_\mu}S^{m_S}$.
\end{lemma}

\proof By definition of the coefficients $\gamma_M$, we have
$$1=\sum_{[M]}[M]\cdot\sum_{[N]}\gamma_N[N]=\sum_{[X]}(\sum_{U\subset X}\gamma_U)[X],$$
where the inner sum runs over all subrepresentations $U$ of $X$. Thus, we have to prove that the coefficients $\gamma_M$ defined as above fulfill the identity
$$\sum_{U\subset M}\gamma_U=0$$
for all $0\not=M\in{\rm mod}_\mu kQ$. Decompose the ${\rm mod}_\mu kQ$-relative socle of $M$ into a direct sum of stable representations $\oplus_{[S]\in{\cal S}_\mu}S^{m_S}$. Since $\gamma_U=0$ if $U$ is not polystable, the above sum over subrepresentations $U$ of $M$ reduces to a sum over the polystable subrepresentations, all of which embed into the (relative) socle of $M$, and thus are of the form $U\simeq\oplus_{[S]\in{\cal S}_\mu} S^{a_S}$ for tuples $(0\leq a_S\leq m_S)_S$. At this point, we have to borrow a remark from the next section: the endomorphism ring ${\rm End}(S)$ of a stable representation is actually a finite extension field of $k$ (see the remark following Definition \ref{absst}). Looking at isotypical components, the number of subrepresentations of a fixed isomorphism class given by a tuple $(a_S)_S$ as above is therefore equal to the cardinality of the product of Grassmannians
$$\prod_{[S]\in{\cal S}}{\rm Gr}_{a_S}({\rm End}(S)^{m_S}),$$
which is
$$\prod_{[S]\in{\cal S}}\left[{m_S\atop a_S}\right]_{|{\rm End}(S)|}.$$
In this last formula, $\left[{m\atop n}\right]_q:=\frac{[m]_q^!}{[n]_q^![m-n]_q^!}$ denotes the $q$-binomial coefficient, where $[n]_q^!=\prod_{i=1}^n[i]_q^!$ is defined via the quantum numbers $[n]_q:=\frac{q^n-1}{q-1}$.\\[1ex]
Thus, we have
\begin{eqnarray*}
\sum_{U\subset M}\gamma_U&=&\sum_{(a_S\leq m_S)_S}\prod_{S}\left[{m_S\atop a_S}\right]_{|{\rm End}(S)|}(-1)^{a_S}|{\rm End}(S)|^{a_S\choose 2}\\
&=&\prod_{S}\sum_{a=0}^{m_S}\left[{m_S\atop a}\right]_{|{\rm End}_A(S)|}(-1)^{a}|{\rm End}(S)|^{a\choose 2}.
\end{eqnarray*}
Now we use the following identity which is easily verified by induction on $m\in{\bf N}$:
$$\sum_{a=0}^m(-1)^aq^{a\choose 2}\left[{m\atop a}\right]=\left\{\begin{array}{ccc}0&,&m\geq 1\\ 1&,&m=0\end{array}\right.$$
Thus, we have $\sum_{U\subset M}\gamma_U=0$ if some $m_S\not=0$, which is the case if and only if $M\not=0$. \hb

\section{Arithmetics of stable representations}\label{arithmetics}

We continue to denote by $k={\bf F}_q$ a fixed finite field. Given a finite field extension $k\subset K$ with Galois group $\Gamma$, we have a corresponding inclusion $kQ\subset KQ$ of path algebras, and an action of $\Gamma$ on $KQ$ by $k$-algebra automorphisms. Thus, we have a base change functor $K\otimes_k\_:{\rm mod} kQ\rightarrow {{\rm mod}} KQ$ and a functorial action of $\Gamma$ on ${\rm mod} KQ$. Given a representation $M$ of $KQ$ and an element $\sigma\in\Gamma$, the representation $M^\sigma$ is defined via the multiplication $a*m=\sigma(a)m$ for all $a\in KQ$ and all $m\in M$. Note, in particular, that $\dimv M^\sigma=\dimv M$, thus $\mu(M^\sigma)=\mu(M)$ for $M\not=0$.\\[1ex]
A representation $M$ of $KQ$ is said to be defined over $k$ if there exists a representation $\overline{M}$ of $kQ$ such that $M\simeq K\otimes_k\overline{M}$. The following is standard:

\begin{lemma}\label{descent} $M$ is defined over $k$ if and only if $M^\sigma\simeq M$ for all $\sigma\in\Gamma$.
\end{lemma}

This motivates the following definition: if $M^\sigma\not\simeq M$ for all $\sigma\in\Gamma$, then $K$ is called a minimal field of definition for $M$. Note that this notion depends on the field $k$ which is fixed throughout.

\begin{lemma}\label{l52} Let $M$ be a representation of $kQ$, and let $X$ be a representation of $KQ$.
\begin{enumerate}
\item If $X$ is stable (resp.~semistable, polystable), then $X^\sigma$ is so, too.
\item $M$ is semistable (resp.~polystable) if and only if $K\otimes_kM$ is so.
\end{enumerate}
\end{lemma}

\proof For the first part, note that a subrepresentation $U\subset X^\sigma$ induces a subrepresentation $U^{\sigma^{-1}}\subset X$. Since dimension types and slopes of representations are preserved under the action of $\Gamma$, all properties follow. For the second part, suppose first that $K\otimes_kM$ is semistable (resp.~polystable). As above, a subrepresentation $U\subset M$ induces a subrepresentation $K\otimes_kU\subset K\otimes_kX$ of the same dimension type (and slope), and the statement follows. To prove the converse, suppose that $M$ is semistable, and consider the individual steps $X_i$ in the Harder-Narasimhan filtration of $K\otimes_kM$. By the first part of the lemma and the defining properties of this filtration, each $X_i$ is stable under the action of $\Gamma$. By Lemma \ref{descent}, there exists a filtration $\overline{X}_i$ of $M$ such that $K\otimes_k\overline{X}_i\simeq X_i$ for all $i$. It is easily seen that the $\overline{X}_i$ define a Harder-Narasimhan filtration of $M$, which is thus trivial by semistability of $M$. Thus, the original filtration by $X_i$ of $K\otimes_kM$ is trivial, too, proving that $K\otimes_kM$ is semistable. Now suppose that $M$ is polystable. We just proved that $K\otimes_kM$ is again semistable, thus we can consider its relative socle $U$. Similar to the above, the polystable representation $U$ is fixed by $\Gamma$, thus descends to a polystable representation $\overline{U}$, which neccessarily equals the relative socle of $M$. But since $M$ is polystable, we have $\overline{U}=M$, thus $U=K\otimes_kM$.\hb
 
\begin{definition}\label{absst} A stable representation $M$ in ${\rm mod}_\mu kQ$ is called absolutely stable if $K\otimes_kM$ is stable for all finite extensions $k\subset K$. An extension field $K\supset k$ is called a splitting field for $M$ if $K\otimes_kM$ is isomorphic to a direct sum of absolutely stable representations.
\end{definition}

\remark Note that, in the second part of the definition, the representation $K\otimes_kM$
is polystable by Lemma \ref{l52}. Note furthermore that the endomorphism ring ${\rm End}(M)$ of a stable representation $M$ is a finite dimensional $k$-division algebra, thus a finite extension field of $k$.\\[1ex]
The following result is proved in the same way as the analogous result for (absolutely) indecomposable representation in \cite{Kac,KR}:

\begin{proposition}\label{psplitting} Let $M$ be a stable representation of $kQ$.
\begin{enumerate}
\item An extension field $K\supset k$ is a splitting field for $M$ if and only if $K$ contains ${\rm End}(M)$.
\item In case $K={\rm End}(M)$, there exists an absolutely stable representation $X$ of $KQ$ with minimal field of definition $K$ such that $K\otimes_kM\simeq\bigoplus_{\sigma\in\Gamma}X^\sigma$.
\end{enumerate}
\end{proposition}

\remark Note that, in the second part of this statement, we have $\dimv M=|\Gamma|\cdot\dimv X$.\\[2ex]
We introduce the following numbers:
\begin{itemize}
\item $s_{d,r}^\mu(q)$ equals the number of isomorphism classes of stable representations $M$ of ${\bf F}_qQ$ of dimension type $d$ with $\dim_{{\bf F}_q}{\rm End}(M)=r$,
\item $a_d^\mu(q)$ equals the number of isomorphism classes of absolutely stable representations of ${\bf F}_qQ$ of dimension type $d$,
\item $\dot{a}_d^\mu(q)$ equals the number of isomorphism classes of absolutely stable representations of ${\bf F}_qQ$ of dimension type $d$ with minimal field of definition ${\bf F}_q$.
\end{itemize}

\begin{proposition}\label{moebius} The following identities hold between the functions just defined:
\begin{enumerate}
\item $s_{d,r}^\mu(q)=\frac{1}{r}\dot{a}_{d/r}^\mu(q^r)$,
\item $a_d^\mu(q^r)=\sum_{s|r}\dot{a}_d^\mu(q^s)$,
\item $s_{d,r}^\mu(q)=\frac{1}{r}\sum_{s|r}\mu(\frac{r}{s})a_{d/r}^\mu(q^s)$.
\end{enumerate}
\end{proposition}

\proof The first identity is just a reformulation of the second part of Proposition \ref{psplitting}, noting that the Galois group $\Gamma$ of an extension $k\subset K$ of degree $r$ is of cardinality $r$. The second identity follows from the definitions. The third identity follows from Moebius inversion applied to the second one, combined with the first one.\hb

\section{The key calculation}\label{key}

After the preparations of the previous section, we can now continue the strategy outlined  at the beginning of section \ref{identity} and apply the integral to the identity of Lemma \ref{keyidentity}.

\begin{proposition} The following identity holds in ${\bf Q}_q[[I]]$:
\begin{samepage} 
$$(\sum_{d\in{\bf N}I_\mu}\frac{|R_d^{ss}({\bf F}_q)|}{|G_d({\bf F}_q)|}t^d)^{-1}=\sum_\xi\prod_{d,r}{{s_{d,r}^\mu(q)}\choose {\sum_m\xi(d,r,m)}}\frac{(\sum_m\xi(d,r,m))!}{\prod_m(\xi(d,r,m)!)}\times$$
\begin{equation}\times\prod_{r,m}\prod_{i=1}^m(1-q^{ri})^{-\sum_d\xi(d,r,m)}t^{\sum_d(\sum_{r,m}m\xi(d,r,m))d},\end{equation}\label{eqN}
\end{samepage}
where the sum on the right hand side runs over all functions $\xi:{\bf N}I_\mu\times{\bf N}_+\times{\bf N}_+\rightarrow{\bf N}$ with finite support.
\end{proposition}

\proof Integration of the identity $e_\mu^{-1}=\sum_{[M]\in{\rm mod}_\mu{\bf F}_qQ}\gamma_M[M]$ of Lemma \ref{keyidentity} yields 
\begin{equation}\int e_\mu^{-1}=\sum_{[M]\in{\rm mod}_\mu{\bf F}_qQ}\frac{\gamma_M}{|{\rm Aut}(M)|}t^{\dimv M}.\end{equation}\label{eq*}
We start by computing the left hand side of this identity:
$$\int e_\mu=\sum_{[X]\, :\, X\in{\rm mod}_\mu {\bf F}_qQ}\frac{1}{|{\rm Aut}(X)|}t^{\dimv X}=\sum_{d\in{\bf N}I_\mu}\frac{|R_d^{ss}({\bf F}_q)|}{|G_d({\bf F}_q)|}t^d.$$
Using Lemma \ref{propint}, this yields
$$\int e_\mu^{-1}=(\int e_\mu)^{-1}=(\sum_{d\in{\bf N}I_\mu}\frac{|R_d^{ss}({\bf F}_q)|}{|G_d({\bf F}_q)|}t^d)^{-1}.$$
We now turn to the right hand side of identity (2). Considering isotypical components of a polystable representation $M=\bigoplus_{[S]\in{\cal S}_\mu}S^{m_S}$, we have
$${\rm Aut}(\bigoplus_{[S]\in{\cal S}_\mu}S^{m_S})\simeq\prod_{[S]}{\rm GL}_{m_S}({\rm End}(S)).$$
Using the identity
$$(-1)^mq^{m\choose 2}|{\rm GL}_m({\bf F}_q)|^{-1}=\prod_{i=1}^m(1-q^i)^{-1},$$
we can thus rewrite the right hand side of (2) as
\begin{equation}\sum_{(m_S)_{[S]}}\prod_{[S]}\prod_{i=1}^{m_S}(1-|{\rm End}(S)|^i)^{-1}t^{\sum_{[S]}m_S\dimv S},\end{equation}\label{eq**}
where the sum runs over all tuples $(m_S)_{[S]\in{\cal S}_\mu}$ with finite support.\\[1ex]
The individual summands in (3) involve only the following data on a given stable representation $S$: its dimension type, the cardinality of its endomorphism ring, and its multiplicity $m_S$ in $M$. This motivates the following: to a polystable representation $M$ of ${\bf F}_qQ$, we associate a function $\xi:{\bf N}I_\mu\times{\bf N}_+\times{\bf N}_+\rightarrow{\bf N}$ with finite support, called its type function, as follows:\\[1ex]
$M$ is of type $\xi$ if it can be written in the form
$$M\simeq\bigoplus_{d,r,m}(S_{d,r,m,1}^m\oplus\ldots\oplus S_{d,r,m,\xi(d,r,m)}^m)$$
for representations $S_{d,r,m,i}$ such that:
\begin{itemize}
\item each $S_{d,r,m,i}$ is stable of dimension type $d$,
\item $|{\rm End}(S_{d,r,m,i})|=q^r$,
\item for each fixed pair $(d,r)$, the $S_{d,r,m,i}$ are pairwise non-isomorphic.
\end{itemize}
We now compute the number $N(\xi)$ of isomorphism classes of polystable representations of type $\xi$:\\[1ex]
Fixing a pair $(d,r)\in{\bf N}I_\mu\times{\bf N}_+$, we have to choose $\sum_m\xi(d,r,m)$ pairwise distinct isomorphism classes $[S_{d,r,m,i}]$ out of $s_{d,r}^\mu(q)$ possible ones (in the notation of the previous section), giving ${{s_{d,r}^\mu(q)}\choose {\sum_m\xi(d,r,m)}}$ choices. For each fixed $m$, the ordering of the $S_{d,r,m,i}$ has to be disregarded, so that we have $\frac{(\sum_m\xi(d,r,m))!}{\prod_m(\xi(d,r,m)!)}$ different choices for the $S_{d,r,m,i}$ out of the given set of $\sum_m\xi(d,r,m)$ isomorphism classes. We thus arrive at the formula
$$N(\xi)=\prod_{d,r}{{s_{d,r}^\mu(q)}\choose {\sum_m\xi(d,r,m)}}\frac{(\sum_m\xi(d,r,m))!}{\prod_m(\xi(d,r,m)!)}.$$
By the definitions, a summand of (3) corresponding to a polystable representation of type $\xi$ reads
$$\prod_{r,m}\prod_{i=1}^m(1-q^{ri})^{-\sum_d\xi(d,r,m)}t^{\sum_d(\sum_{r,m}m\xi(d,r,m))d}.$$
Thus, we arrive at the claimed identity.\hb\\[3ex]
Next we derive a recursive formula for the numbers $a_d^\mu(q)$ of absolutely stable representations defined in the previous section from the above identity (1).\\[1ex]
First of all, by Proposition \ref{moebius}, we know that $s_{d,r}^\mu(q)=0$ unless $r|d$, in which case $s_{d,r}^\mu(q)=\frac{1}{r}\sum_{s|r}\mu(\frac{r}{s})a_{d/r}^\mu(q^s)$. We can rewrite identity (1) accordingly.\\[1ex]
We define the dimension type of a type function $\xi$ by $$\dimv \xi:=\sum_d(\sum_{r,m}m\xi(d,r,m))d\in{\bf N}I_\mu.$$

Given $(d,r,m)\in{\bf N}I_\mu\times{\bf N}_+\times{\bf N}_+$, denote by $\xi_{d,r,m}$ the function taking value $1$ on $(d,r,m)$ and $0$ otherwise. The summand on the right hand side of the above identity (1) corresponding to the function $\xi_{d,1,1}$ is easily seen to equal $(1-q)^{-1}a_d^\mu(q)t^d$.\\[1ex]
The following is straightforward:

\begin{lemma}\label{trivrec} Let $a=\sum_{d\in{\bf N}I}\alpha_dt^d$ be an element of ${\bf Q}_q[[I]]$ such that $\alpha_0=1$. Writing $a^{-1}=\sum_d\beta_dt^d$, the $\beta_d\in{\bf Q}$ are given recursively by $\beta_0=1$ and, for all $d\in{\bf N}I\setminus\{0\}$:
$$\beta_d=-\sum_{e\lneqq d}q^{\langle e-d,e\rangle}\alpha_{d-e}\beta_e.$$
\end{lemma}

Applying this lemma to the left hand side of identity (1) and using the above discussion, we see that, for each $d_0\in{\bf N}I$, we have
\begin{equation} a_{d_0}^\mu(q)=(1-q)\cdot(\beta_{d_0}-\sum_{{\dimv\xi=d_0,}\atop{\xi\not=\xi_{d_0,1,1}}}\prod_{r|d}{{\frac{1}{r}\sum_{s|r}\mu(\frac{r}{s})a_{d/r}^\mu(q^s)} \choose {\sum_m\xi(d,r,m)}} \ldots).\end{equation}\label{eq+}

\begin{lemma} Only terms containing the $a_e^\mu(q^s)$ for $e\lneqq d_0$ contribute to the sum in identity (4).
\end{lemma}

\proof If a term $a_{d/r}^\mu(q^s)$ contributes to the summand corresponding to a type function $\xi$, we have $c:=\sum_m\xi(d,r,m)\geq 1$, so that $\sum_{r,m}m\xi(d,r,m)\geq c$. This implies
$$d_0=\dimv\xi=\sum_d(\sum_{r,m}m\xi(d,r,m))d\geq cd\geq d.$$
Therefore, $d/r\leq d\leq d_0$, and $d/r=d_0$ only if $r=1$ and $c=1$. In this latter case, it is easily seen that $\xi=\xi_{d,1,1}$, a contradiction.\hb

For computation of the $a_d^\mu(q)$, it thus remains to compute the $\beta_d$. This is contained in \cite{HNS}, namely Corollary 6.2 there (note that the assumption of \cite{HNS} that the quiver $Q$ does not have oriented cycles is not essential for this result):\\[1ex]
For all dimension types $d\in{\bf N}I$, we have
$$\frac{|R_d^{ss}({\bf F}_q)|}{|G_d({\bf F}_q)|}=\sum_{d^*}(-1)^{s-1}q^{-\sum_{i<j}\langle d^j,d^i\rangle}\prod_{i=1}^s\frac{|R_{d^i}({\bf F}_q)|}{|G_{d^i}({\bf F}_q)|},$$
where the sum runs over all tuples $d^*=(d^1,\ldots,d^s)$ of dimension types such that $d=\sum_id^i$ and $\mu(\sum_{j=1}^id^j)>\mu(d)$ for all $i<s$.\\[2ex]
Now the varieties $R_d({\bf F}_q)$ and $G_d({\bf F}_q)$) are just an affine space (resp.~a product of general linear groups), so their cardinalities are known. This allows calculation of all $\beta_d$ by the recursion Lemma \ref{trivrec}, and thus of the $a_d^\mu(q)$ by the above discussion. We summarize the above discussion:

\begin{theorem}\label{doublerecursion} With the above notation, the number $a_d^\mu(q)$ of absolutely stable representations of ${\bf F}_qQ$ of dimension type $d$ is given by the recursive formula
$$a_d^\mu(q)=(1-q)\cdot(\beta_d-\sum_{{\dimv\xi=d,}\atop{\xi\not=\xi_{d,1,1}}}\prod_{r|d}{{\frac{1}{r}\sum_{s|r}\mu(\frac{r}{s})a_{d/r}^\mu(q^s)}\choose {\sum_m\xi(d,r,m)}}\frac{(\sum_m\xi(d,r,m))!}{\prod_m(\xi(d,r,m)!)}\times$$
$$\times\prod_{r,m}\prod_{i=1}^m(1-q^{ri})^{-\sum_d\xi(d,r,m)}t^{\sum_d(\sum_{r,m}m\xi(d,r,m))d}),$$
where the sum runs over type functions $\xi$, and the $\beta_d$ are given by the recursive formula
$$\beta_d=-\sum_{e\lneqq d}q^{\langle e-d,e\rangle}(\sum_{f^*}(-1)^{s-1}q^{-\sum_{i<j}\langle f^j,f^i\rangle}\prod_{i=1}^s\frac{|R_{f^i}({\bf F}_q)|}{|G_{f^i}({\bf F}_q)|})\cdot\beta_e,$$
the sum running over all tuples $f^*=(f^1,\ldots,f^s)$ of dimension types such that $d-e=\sum_if^i$ and $\mu(\sum_{j=1}^if^j)>\mu(d-e)$ for all $i<s$.
\end{theorem}

From these formulas, we can see immediately:

\begin{corollary}\label{ratfun} For all $d\in{\bf N}I$, there exists a rational function $A_d^\mu(t)\in{\bf Q}(t)$ such that, for all finite fields $k$, we have $A_d^\mu(|k|)=a_d^\mu(|k|)$.
\end{corollary}

\section{Polynomials counting rational points}\label{counting}

As before, let $k$ be a fixed finite field with $q$ elements, and fix an algebraic closure $\overline{k}$ of $k$.

\begin{proposition}\label{poly} Let $X$ be an irreducible $\overline{k}$-variety, and let $X_0$ be a variety over $k$ such that $X\simeq \overline{k}\times_kX_0$. Suppose there exists a rational function $P_X(t)\in{\bf Q}(t)$ such that $P_X(|K|)$ equals the number $|X_0(K)|$ of $K$-rational points of $X_0$ for all finite extensions $K\supset k$. Then:
\begin{enumerate}
\item $P_X(t)\in{\bf Z}[t]$,
\item $P_X(1)$ equals the Euler characteristic $\chi_c(X)$ of $X$ in $\ell$-adic cohomology with compact support.
\end{enumerate}
\end{proposition}

\proof Write $P_X(t)=\frac{1}{N}Q(t)+\frac{R(t)}{S(t)}$ for $N\in{\bf N}$ and polynomials $Q,R,S\in{\bf Z}[t]$ such that $\deg R<\deg S$. Suppose that $P_X$ is a proper rational function, i.e.~$S$ is not constant. Since $\lim_{t\rightarrow\infty}\frac{R(t)}{S(t)}=0$, we can find $r\gg 0$ such that $0<|\frac{R(q^r)}{S(q^r)}|<\frac{1}{N}$. Then $0<|NP_X(q^r)-Q(q^r)|<1$, contradicting $NP(t)-Q(t)\in{\bf Z}[t]$. Thus, $P_X(t)\in{\bf Q}[t]$, which we write as $P_X(t)=\sum_mc_mt^m$.\\[1ex]
For $i=0\ldots 2\dim X$, let $\alpha_{i,1},\ldots,\alpha_{i,b_i}$ be the eigenvalues (listed with mul\-ti\-pli\-cities) of Frobenius acting on the $i$-th $\ell$-adic cohomology with compact support $H_c^i(X,\overline{{\bf Q}_\ell})$ of $X$ for some prime $\ell\not={\rm char} k$. By Grothendieck's trace formula for Frobenius \cite{De} and the assumptions, we have for each extension field $K\supset k$ of degree $r$ the following identity:
$$\sum_mc_mq^{rm}=P_X(|K|)=|X_0(K)|=\sum_{i=0}^{2\dim X}(-1)^i\sum_{j=1}^{b_i}\alpha_{i,j}^r$$
A standard Vandermonde argument yields the following:\\[1ex]
Let $\beta_n,\gamma_n$ for $n=1,\ldots,t$ be complex numbers such that all $\gamma_n$ are non-zero and pairwise distinct. Assume that $\sum_{n=1}^t\beta_n\gamma_n^r=0$ for all $r\geq 1$. Then all $\beta_n$ are zero.\\[1ex]
Applying this to the above identity, with the $q^m$ and the $\alpha_{i,j}$ playing the role of the $\gamma_n$, we derive:
\begin{itemize}
\item each $\alpha_{i,j}$ is of the form $q^m$ for some $m$,
\item $c_m=\sum_i(-1)^i|\{j\, :\, \alpha_{i,j}=q^m\}|\in{\bf Z}$.
\end{itemize}
Thus, $P_X(t)\in{\bf Z}[t]$, and $P_X(1)=\sum_mc_m=\sum_{i}(-1)^i\sum_m|\{j\, :\, \alpha_{i,j}=q^m\}|=\sum_i(-1)^ib_i=\chi_c(X)$.\hb

We want to apply this result to the moduli space $M^s_d(Q)(\overline{k})$. To do this, we first define all geometric objects considered above as schemes over ${\bf Z}$. We define the affine scheme ${R}_d=\bigoplus_{\alpha:i\rightarrow j}{\rm Hom}_{\bf Z}({\bf Z}^{d_i},{\bf Z}^{d_j})$ with the algebraic group scheme ${ G}_d:=\prod_{i\in I}{\rm GL}_{d_i}({\bf Z})$ acting on it as before. Using the stability $\Theta$, we define open subschemes ${ R}_d^{s}\subset{ R}_d^{ss}\subset{ R}_d$ as in \cite[2.1]{CBVdB} (where they are defined using Seshadri's generalization of Geometric Invariant Theory to arbitrary base rings \cite{Ses}). Again using \cite{Ses}, we can define a scheme ${ M}_d^s(Q)$, whose geometric points over the algebraically closed field $\overline{k}$ are in bijection with the orbits of $G_d(\overline{k})$ in $R_d^s(\overline{k})$, that is, to the isomorphism classes of stable representations of $\overline{k}Q$ of dimension type $d$.\\[1ex]
For the finite field $k$, let $M_d^s(Q)_k$ be the $k$-variety obtained by base change to $k$ from ${ M}_d^s(Q)$. From the above, we see that $\overline{k}\times_kM_d^s(Q)_k$ is isomorphic to the moduli space $M_d^s(Q)(\overline{k})$ defined in section \ref{quivermoduli}. From this construction, it is clear that, for any finite extension $K\supset k$, the $K$-rational points of $M_d^s(Q)_k$ correspond precisely to the isomorphism classes of absolutely stable representation of $KQ$ of dimension type $d$.\\[2ex]
Applying Proposition \ref{poly} to the varieties $X=M^s_d(Q)(\overline{k})$ and $X_0=M_d^s(Q)_k$, we get the main result of the paper.

\begin{theorem}\label{polynomiality} The rational function $A_d^\mu(t)$ defined in Corollary \ref{ratfun} is already a polynomial in ${\bf Z}[t]$, whose evaluation at $t=1$ gives the Euler characteristic in $\ell$-adic cohomology with compact support of $M^s_d(Q)(\overline{k})$.
\end{theorem}

\section{Examples and comments}\label{examples}

Although the double recursion Theorem \ref{doublerecursion} is difficult to use for explicit computations by hand, it is easily implemented on a computer (the part of the formula corresponding to the explicit formula \cite[Corollary 6.2]{HNS} should be replaced by the algorithm \cite[Corollary 6.9]{HNS} for this purpose). All explicit formulas below were found with the aid of such an implementation.\\[2ex]
First we consider the case of an arbitrary quiver $Q$ with trivial stability $\Theta=0$. In this case, it follows from the definitions that all representations of $Q$ are semistable of slope $0$, and the stable (resp.~polystable) representations are just the simple (resp.~semisimple) ones. Thus, the polynomial $A_d^0(t)$ counts isomorphism classes of absolutely simple representations of $Q$ of dimension vector $d$ over finite fields.\\[1ex]
This leads to a simplification of the double recursion
Theorem \ref{doublerecursion}. Namely, the left hand side of identity (1) of section \ref{key} reads $(\sum_{d\in{\bf N}I}\frac{|R_d({\bf F}_q)|}{|G_d({\bf F}_q)|}t^d)^{-1}$ in this case. Thus, no reference to the formula of \cite{HNS} has to be made.\\[2ex]
We will work out this recursion in more detail in the case of the quiver $Q_m$ consisting of a single vertex $i$ and $m$ loops. The polynomal $A_d^{(m)}(t):=A_d^0(t)$ for this quiver thus counts absolutely simple representations of a free algebra in $m$ generators. In linear algebra terms, we have:\\[1ex]
For any finite field $k$, the evaluation $A_d^{(m)}(|k|)$ equals the number of simultaneous conjugacy classes of $m$-tuples of $d\times d$-matrices with entries in $k$, whose scalar extension to $\overline{k}$ does not admit a non-trivial common invariant subspace.\\[1ex]
We first note that ${\bf Q}_q[[I]]$ can be identified with the formal power series ring ${\bf Q}[[t]]$ (with the usual multiplication) by identifying $t^d\in{\bf Q}_q[[I]]$ with the element $q^{(1-m){d\choose 2}}t^d\in{\bf Q}[[t]]$ for all $d$. The recursive description of the inverse of a series in ${\bf Q}_q[[I]]$ in Lemma \ref{trivrec} can be replaced by the following explicit formula for inverses in ${\bf Q}[[t]]$, whose proof is again straightforward:
$$(1+\sum_{d=1}^{\infty} a_dt^d)^{-1}=\sum_{\lambda\in\Lambda}(-1)^{l(\lambda)}{{l(\lambda)}\choose{\mu_1(\lambda),\mu_2(\lambda),\ldots}}\prod_{i}a_{\lambda_i}t^{|\lambda|}.$$
In this formula, $\Lambda$ denotes the set of partitions, and for a partition $\lambda$, we have the notation $l(\lambda)$ for its length, $|\lambda|=\sum_i\lambda_i$ for its weight, and $\mu_i(\lambda)$ for the number of occurences of $i$ in $\lambda$.\\[1ex]
Using these two simplifications, we can rewrite (after some calculations) the left hand side of identity (1) as
$$\sum_{\lambda\in\Lambda}(-1)^{l(\lambda)}{{l(\lambda)}\choose{\mu_1(\lambda),\mu_2(\lambda),\ldots}}\prod_{i=1}^{l(\lambda)}\frac{q^{m{{\lambda_i+1}\choose 2}}}{(q^{\lambda_i}-1)\ldots(t-1)}t^{|\lambda|}\in{\bf Q}[[t]].$$
We can simplify the right hand side of identity (1) by resolving a slight redundance: if a type function $\xi$ contributes to the sum, then we already have $r|d$ whenever $\xi(d,r,m)\not=0$, since $s_{d,r}^\mu(q)=0$ otherwise. This suggests replacing $\xi$ by the function $\tilde{\xi}$ given by $\tilde{\xi}(d,r,m):=\xi(rd,r,m)$. The neccessary modifications to the right hand side of identity (1) resulting from this reparametrization are easily performed. After some more calculations, we finally arrive at the following recursion:

\begin{theorem}\label{freealgebra} For all $d\in{\bf N}$, we have $A_d^{(m)}(t)$=
$$(1-t)\left((t^{(m-1){d\choose 2}}\sum_{\lambda\in\Lambda_d}(-1)^{l(\lambda)}{{l(\lambda)}\choose{\mu_1(\lambda),\mu_2(\lambda),\ldots}}\prod_{i=1}^{l(\lambda)}\frac{t^{m{{\lambda_i+1}\choose 2}}}{(t^{\lambda_i}-1)\ldots(t-1)}\right.$$
$$\left. -\sum\prod_{i,j,k}\frac{1}{\xi_{ijk}!}(\prod_{l=1}^k(1-t^{lj}))^{-\xi_{ijk}}\prod_{i,j}{{\frac{1}{j}\sum_{r|j}\mu(\frac{j}{r})a_i^{(m)}(t^r)}\choose{\sum_k\xi_{ijk}}}(\sum_k\xi_{ijk})!\right),$$
the second sum running over all functions $\xi:{\bf N}_+^3\rightarrow{\bf N}$ such that $\sum_{i,j,k}ijk\xi_{ijk}=d$ and $\xi\not=\delta_{d,1,1}$, and $\Lambda_d$ denoting the set of partitions of $d$.
\end{theorem}

We give examples in low dimensions (the formula for dimension $4$ is already rather unpleasant):\\[1ex]
\example We have $A_1^{(m)}(t)=t^m$, $A_2^{(m)}(t)=t^{2m}\frac{(t^m-1)(t^{m-1}-1)}{t^2-1}$ and $A_3^{(m)}(t)=$
$$t^{3m+1}\frac{(t^m-1)(t^{2m-2}-1)(t^{3m-2}+t^{2m-2}-t^m-2t^{m-1}-t^{m-2}+t+1)}{(t^3-1)(t^2-1)}.$$

As another prominent example of quiver moduli, we consider the $m$-arrow Kronecker quiver $Q=i\stackrel{(m)}{\Rightarrow}j$ with dimension type $d=i+j$ and stability $\Theta=i^*$. Then we can find
$$A_d^\mu(t)=\frac{t^m-1}{t-1}\mbox{ and }A_{2d}^\mu(t)=\frac{t^2(t^m-1)(t^{m-1}-1)^2(t^{m-2}-1)}{(t^2-1)^2(t-1)}.$$

Finally, we consider one of the classical examples of Geometric Invariant Theory, namely the space of stable configurations of $m$ points in the projective line modulo projective transformations. As in \cite[2. Example B]{HNS}, we consider the quiver with vertices $i_0,i_1,\ldots,i_m$ and arrows $i_k\rightarrow i_0$ for all $k$. The dimension type is $d=2i_0+i_1+\ldots+i_m$, the stability is $\Theta=-i_0^*$. Generalizing results of \cite{M} for the case of odd $m$ (where all semistables are already stable), we get
$$A_d^\mu(t)=\left\{\begin{array}{lll}\frac{1}{q(q-1)}((q+1)^{2k}-\sum_{l=0}^k{2k+1\choose l}q^l)&,&m=2k+1\\
\frac{1}{q(q-1)}((q+1)^{2k-1}-\sum_{l=0}^k{2k\choose l}q^l+\frac{1}{2}{2k\choose k}q)&,&m=2k\end{array}\right.$$

One of the main applications of Theorem \ref{doublerecursion} is the opportunity for experiments in low-dimensional cases. As a typical example, computer experiments with the formula Theorem \ref{freealgebra} for dimensions $d=2,\ldots,8$ led to the following prediction: the polynomial $A_d^{(m)}(t)$ has a zero of order $1$ at $t=1$, and $$\frac{A_d^{(m)}(t)}{t-1}|_{t=1}=\frac{1}{d}\sum_{r|d}\mu(\frac{d}{r})m^r.$$
In view of Proposition \ref{poly}, this should be properly formulated as a statement about the Euler characteristic of an appropriate $k^*$-quotient of ${M}_d^0(Q_m)(k)$. Such a statement will be established by entirely different (localization) techniques in \cite{LQM}.

\end{document}